\documentclass{amsart}

\usepackage{epsfig}
\usepackage{amsthm}
\usepackage{amssymb}
\usepackage{amsmath}
\usepackage{amscd}

\newcommand{\labell}[1] {\label{#1}}
\newtheorem {Theorem}   {Theorem} 

\numberwithin{Theorem}{section}

\newtheorem {Lemma}[Theorem]    {Lemma}         
\newtheorem {Proposition}[Theorem]{Proposition}  
\theoremstyle{definition}

\theoremstyle{remark}
\newtheorem{Remark}[Theorem]{Remark}

\newtheorem{Claim}[Theorem]{Claim}
\expandafter\chardef\csname pre amssym.def at\endcsname=\the\catcode`\@ 
\catcode`\@=11 
\def\undefine#1{\let#1\undefined} 
\def\newsymbol#1#2#3#4#5{\let\next@\relax 
 \ifnum#2=\@ne\let\next@\msafam@\else 
 \ifnum#2=\tw@\let\next@\msbfam@\fi\fi 
 \mathchardef#1="#3\next@#4#5}
\def\mathhexbox@#1#2#3{\relax 
 \ifmmode\mathpalette{}{\m@th\mathchar"#1#2#3}% 
 \else\leavevmode\hbox{$\m@th\mathchar"#1#2#3$}\fi} 
\def\hexnumber@#1{\ifcase#1 0\or 1\or 2\or 3\or 4\or 5\or 6\or 7\or 8\or 
 9\or A\or B\or C\or D\or E\or F\fi} 
 
\font\teneufm=eufm10 
\font\seveneufm=eufm7
\font\fiveeufm=eufm5 
\newfam\eufmfam
\textfont\eufmfam=\teneufm 
\scriptfont\eufmfam=\seveneufm 
\scriptscriptfont\eufmfam=\fiveeufm 
\def\frak#1{{\fam\eufmfam\relax#1}} 
\catcode`\@=\csname pre amssym.def at\endcsname 

\def	\eps	{\epsilon}

\def	\C      {{\mathbb C}}
\def	\R	{{\mathbb R}}
\def	\Z	{{\mathbb Z}}
\def	\N	{{\mathbb N}}
\def	\Q	{{\mathbb Q}}

\def	\T	{{\mathbb T}}

\def    \ra     {{\rightarrow}}
\def    \haf    {{\frac{1}{2}}}

\def	\codim	{\operatorname{codim}}

\def	\SB	{\operatorname{SB}}
\def	\CL	{\operatorname{CL}}
\def    \l      {\langle}
\def    \r      {\rangle}

\begin{document}

%%%%%%%%%%%%%%%%%%%%%%%%%%%%%%
%   TEXT FORMATTING

%%%%%%%%%%%%%%%%%%%%%%%%%%

%%%%%%%%%%%%%%%%%%%%%%%%%%

%%%%%%%%%%%           BEGINNING OF  TEXT

%%%%%%%%%%%%%%%%%%%%%%%%%%

\title[Periodic orbits near symplectic extrema]{Periodic orbits of Hamiltonian flows near symplectic extrema }

\author[Viktor L. Ginzburg]{Viktor L. Ginzburg}
\author[Ely Kerman]{Ely Kerman}
\address{Department of Mathematics, UC Santa Cruz, 
Santa Cruz, CA 95064, USA; The Fields Institute, 222 College Street, Toronto, Ontario M5T 3J1, Canada}
\email{ginzburg@math.ucsc.edu; ekerman@fields.utoronto.ca}

\date{\today}

\thanks{The work is partially supported by the NSF and by the faculty
research funds of the University of California, Santa Cruz.}


\bigskip

\begin{abstract}
For Hamiltonian flows we establish the existence of periodic orbits on a 
sequence of level sets approaching a Bott-nondegenerate symplectic extremum 
of the Hamiltonian. As a consequence, we show that a charge on a compact
manifold with a nondegenerate (i.e. symplectic) magnetic field has periodic
orbits on a sequence of energy levels converging to zero.
\end{abstract}

\maketitle

\section{Introduction} 

In the early seventies, Alan Weinstein proved the following result which was subsequently reproved by Jurgen Moser (using different methods) and is now known as the Weinstein-Moser Theorem (see \cite{we1,we3,mo2}).

\begin{Theorem}[Weinstein-Moser]
\labell{thm:wm}
Let $H$ be a smooth function on a symplectic manifold of dimension $2n$. Then 
the Hamiltonian flow of $H$ has at least $n$ periodic orbits on all level sets sufficiently close to a nondegenerate extremum point of $H$.
\end{Theorem}
In this paper, we attempt to extend this result from extremal points to higher dimensional extrema of the Hamiltonian. In particular, we consider the case of symplectic extremal submanifolds. More precisely, let $H$ be a smooth function on a symplectic manifold $(W,\Omega)$ such that $H$ reaches an extremum at a compact symplectic Bott-nondegenerate submanifold $M^{2l} \subset W^{2n}.$ We prove the following.   
\begin{Theorem}
\labell{thm:main}
The Hamiltonian flow defined on $(W,\Omega)$ by the function $H$ has at least one periodic orbit on a sequence of  energy levels converging to $M$.
\end{Theorem}

In comparing these theorems, we see that here the existence of periodic orbits for all sufficiently close level sets is weakened to existence on a sequence of level sets approaching $M$. We also note that the lower bound for the number of periodic orbits is replaced by simple existence. 
It is unlikely that this result is sharp. (The uncertainty as to where the true boundaries of such existence results should lie is indicative of a lack of examples in this area.) However, other existence results of this kind have recently been obtained, see \cite{po} (Theorem \ref{thm:po} below) and \cite{le}. Moreover, a result that is similar in strength to the Weinstein-Moser Theorem can be proved if one imposes certain compatibility conditions on $\Omega$ and the Hessian of $H$ on $M$. Specifically, in \cite{ke} it is shown that under such assumptions there are at least $\CL(M,\Q)+(n-l)$ periodic orbits on all level sets sufficiently close to $M$. Here $\CL(M,\Q)$ denotes the cup-length of
$M$ over $\Q$.

\subsection{Symplectic magnetic flows}
The question addressed here is motivated further by the following interesting set of examples. Let $(M,\omega)$ be a compact symplectic manifold and $g$ a Riemannian metric on $M$. Consider the Hamiltonian flow defined on $T^*M$ by the kinetic energy Hamiltonian 
\begin{align*}
H_g \colon T^*M &\to \R \\
(q,p) &\mapsto \|p\|_{g^{-1}}^2
\end{align*}  
and the twisted symplectic from $d\lambda + \pi^*\omega$. Here $\pi \colon T^*M \to M$ is the projection map and $\lambda$ is the canonical Liouville one-form. These flows describe the  motion of a charged particle in a nondegenerate magnetic field and will be referred to as symplectic magnetic flows. The zero section of $T^*M$ is a symplectic minimum of $H_g$ and in the context of the result above we are concerned with the existence of periodic orbits on low energy levels. For such flows, Theorem \ref{thm:main} implies the following result.

\begin{Theorem}
For any symplectic form $\omega$ and metric $g$ on $M$, the corresponding symplectic magnetic flow has periodic orbits on a sequence of low energy levels converging to zero.
\end{Theorem}

To the knowledge of the authors, this is the most general existence result for symplectic magnetic flows. However, much stronger results have been established in a variety of different cases.

\begin{Theorem}[\cite{ar1,gi1}]
Let $M$ be a surface of genus $k$. Then for any choice of $\omega$ and $g$, the corresponding symplectic magnetic flow has at least three periodic orbits on all sufficiently low energy levels and at least $2k+2$ if they are nondegenerate.
\end{Theorem}

 \begin{Theorem}[\cite{ke, gk1}]
Let the metric $g$ be of the form $\omega(\cdot\,, \,J\cdot)$ for some almost complex structure $J$ on $M$. Then the corresponding symplectic magnetic flow has at least $\CL(M,\R)+l$ periodic orbits on all sufficiently low energy levels and at least $\SB(M)$ if all orbits are nondegenerate. Here $\SB(M)$ denotes the sum of Betti numbers of $M$.
\end{Theorem}

In fact, the construction above produces a Hamiltonian flow from any manifold $M$, closed two-form $\omega$ and metric $g$, i.e. $\omega$ can be degenerate. We call these flows magnetic flows and refer the reader to \cite{gi3} for a more detailed discussion of them and further references. The following are recent results for magnetic flows which also hold in the symplectic case.

\begin{Theorem}[\cite{gk1}]
For any closed two-form $\omega$ and metric $g$ on $\T^n$, the corresponding magnetic flow has periodic orbits on almost all energy levels. (In fact, bounded neighborhoods of the zero section have finite Hofer--Zehnder capacity.) 
\end{Theorem}

\begin{Theorem}[\cite{le,po}]
\labell{thm:po}
For any metric on $M$ and any nonzero weakly-exact\footnote{Recall that a
form $\omega$ is called weakly-exact if $\omega$ is closed and
$[\omega]|_{\pi_2(M)}=0$.} two-form, the 
corresponding magnetic flow has contractible periodic orbits on a sequence 
of energy levels converging to zero.
\end{Theorem}

These last two theorems represent interesting applications of tools from other areas of symplectic topology to the existence question. The first result is a direct application of the work on the Hofer--Zehnder capacity in \cite{fhv} and \cite{ji}. The second result is proved  for $M=\T^n$ by Polterovich using Hofer's metric on the space of Hamiltonian diffeomorphisms (see \cite{po}). In particular, he utilizes the relation between the nonminimizing geodesics of Hofer's metric and the existence of contractible periodic orbits. This work is then extended to the more general form, as stated here, by Macarini in \cite{le}.
 
\section{Limiting dynamics and the variational problem} 

\subsection{The variational problem}
Before starting the proof of Theorem \ref{thm:main}, we recall the variational framework for proving the existence of 
periodic orbits of the Hamiltonian flow of a function $H$
on a symplectic manifold $(W, \Omega)$. First we choose a suitable class of loops in $W$, say the Fr\'{e}chet manifold of smooth loops, $C^{\infty}(S^1,W)$. Then on $C^{\infty}(S^1,W)$ we consider the one-form $\frak{F}$ which takes $v$, an element of the tangent space at $\sigma \in C^{\infty}(S^1,W)$, to
\begin{equation}
\label{eq:frakF}
 \frak{F}(v) =  \int_0^{1} \Omega(\sigma(s))(\dot{\sigma}(s), v(s)) \, ds - \int_0^{1} dH({\sigma}(s))(v(s))\,ds.
\end{equation}
It is clear from the least action principle that the zeroes of $\frak{F}$ on $C^{\infty}(S^1,W)$ are periodic orbits of the Hamiltonian flow with period equal to one.

At this point we are faced with several difficulties. First among them is the fact that the variational problem is far more tractable if $\frak{F}$ is exact. For then we may look for the critical points of a functional instead of the zeros of a one-form. As well, there need not exist closed orbits with period equal to one. Finally, as stated, the variational approach does not allow us to search for periodic orbits on a fixed level set of $H$.

The analysis of the limiting dynamics given below suggests that the search for periodic orbits on level sets near $M$ can be restricted to the subset of small loops in $C^{\infty}(S^1,W)$ that lie close to $M$. In making this restriction, we are able to overcome the first problem by finding an action functional whose derivative is equal to $\frak{F}$ on this subset. The restriction also  allows us to considerably simplify the analytic setting by using the method developed by Weinstein in \cite{we4} to prove the Arnold conjecture for $C^0$-small Hamiltonians. The other difficulties are surmounted by adapting a set of techniques developed by Viterbo, and Hofer and Zehnder in 
\cite{hz1,hz2,hz3,vi}. Namely, we modify the
function $H$ (and hence our functional) to force any periodic orbits on the desired energy levels to have positive action. Then
we show that the new action functional has
critical points with positive action by using a ``linking argument'' as in \cite{hz1,hz2,hz3}.

\subsection{The limiting dynamics}

For the case when $M$ is an extremal point of the function $H \colon W^{2n} \ra \R,$ we recall briefly how one begins to look for periodic orbits of the 
Hamiltonian flow of $H$ on level sets near $M$. Darboux's theorem allows us to work in a neighborhood of the origin in $\R^{2n}$ with the canonical symplectic form $\Omega_0$. Here the function looks like 
$$ 
H(z)= \haf H_{zz}(0)z^2 + O(z^3)
$$ 
and we may assume that the origin is a local minimum so that the quadratic form  $H_{zz}(0)$ is positive definite. Setting 
\begin{align*}
H_{\eps}(z) &=\eps^{-2}H(\eps z) \\
            & = \haf H_{zz}(0)z^2 + O(\eps)
\end{align*}
it is easily checked that for $\eps\neq 0$ the flow of $X_H$ on $\{H= \eps^2\}$ is diffeomorphic to the flow of the Hamiltonian vector field $X_{H_{\eps}}$ on $\{H_{\eps}= 1\},$ via the rescaling map $z \mapsto \eps z$. 

The vector field $X_{H_{\eps}}$ is defined by 
$$
i_{X_{H_{\eps}}} \Omega_0 = d {H_{\eps}}
$$ 
and we see that $X_{H_\eps}$ is a Hamiltonian perturbation of the linear Hamiltonian vector field $X_0$ given by the equation 
$$
i_{X_{0}} \Omega_0 = H_{zz}(0)z.
$$ 
Since $H_{zz}(0)$ is positive definite,
there exists a change of variable $z \mapsto y$ which preserves $\Omega_0$ and puts $H_{zz}(0)$ in the form 
$$ 
H_{zz}(0)y^2 =\sum_{i=1}^{n} a_i(y_i^2+y_{i+n}^2),
$$ 
see \cite[\S 1.7]{hz3}.
In these coordinates it becomes clear that $X_0$ describes the quasiperiodic motion of $n$ uncoupled harmonic oscillators with frequencies $a_i$. The problem then reduces to showing that the periodic ``normal modes'' of the linear system persist under the perturbation $X_{H_{\eps}}$ (see \cite{ly,mo2,we1,we3} and \cite{fara}).

In replacing a critical point by a critical submanifold $M$ we will rescale globally in the normal directions to $M$ and show that we still get a well-defined and useful limiting vector field.

Setting $H|_M=0$, the Tubular Neighborhood Theorem allows us to assume that, 
for sufficiently small $\epsilon>0$, the level sets $\{H=\epsilon^2\}$ lie in 
a neighborhood of the zero section in the total space of a normal bundle $N$
to $M$. Accordingly, we may replace the manifold $W$ by this neighborhood which
we will still denote by $W$.
We choose the normal bundle to be $(TM)^{\Omega}$, the symplectic orthogonal complement to $TM$. By Weinstein's Symplectic Neighborhood Theorem (see \cite{we2}) we may also assume that $\Omega$ restricts to the fibres in $(TM)^{\Omega} \cap W$ as a constant linear symplectic form, $\Omega^N$. Thus, the level sets of interest lie in a symplectic vector bundle which is also equipped with a fibrewise positive-definite quadratic form $d^2_N H$, given by the Hessian of $H$ in the normal directions to $M$. As above, there exist coordinates
$\{y_i(x)\}_{i=1}^{2(n-l)}$ in each fibre $E_x$ such that $\Omega^N (x)$ is the canonical symplectic form on ${\R}^{2(n-l)}$ and
\begin{equation}
\labell{eq:quasi}
 d^2_N H(x)y^2= \sum_{i=1}^{(n-l)} a_i(x)(y_i^2+y_{i+(n-l)}^2).
\end{equation}
Note that in general these coordinates are not unique and cannot be chosen 
to depend smoothly or even continuously on $x$. However, the eigenvalues 
$a_i(x)$ of $d^2_N H (x)$ with respect to $\Omega^N(x)$ are well defined.
 
Starting with the Hamiltonian dynamical system defined on $W \subset N$ 
by 
$$
i_{X_H}\Omega=dH
$$
we let $\Phi\colon N \ra N$ be the global fibrewise dilation by a factor of 
$\epsilon$, and set $$X_{\eps}:={\Phi^{-1}}_* X_H$$ and $$\tilde{\Omega}_{\eps} :=\epsilon^{-2} \Phi^*\Omega.$$ After this rescaling the new Hamiltonian dynamical system is given by 
\begin{equation}
\labell{eq:dyn}
i_{X_{\eps}}\tilde{\Omega}_{\eps}=d(\epsilon^{-2}\Phi^* H).
\end{equation} 

\begin{Lemma}
\labell{Lemma:lemma1}

As $\epsilon \ra 0$, $X_{\eps}$ approaches a fibrewise, quasiperiodic vector field $X_0$. In particular, in each fibre $E_x$, the vector field $X_0$ is the (linear) Hamiltonian vector field of the positive-definite form $d^2_N H(x)$ with respect to the symplectic form ${\Omega}^N(x)$. 

\end{Lemma}

\begin{proof}

One can check that the fibre components of the $\tilde{\Omega}_{\eps}$ are independent of $\eps$ and equal to ${\Omega}^N.$ However, the limit of the $\tilde{\Omega}_{\eps}$ as $\eps \ra 0$ is not defined, i.e. in coordinates, all the terms with components along the base blow up in the limit. Using the bundle isomorphisms 
$$ 
(\tilde{\Omega}_{\eps})^{\flat} \colon TW \ra T^*W,
$$ 
defined by each of the nondegenerate forms $\tilde{\Omega}_{\eps}$, we can construct the dual bivectors $-\tilde{\Omega}_{\eps}^{-1} \in \Lambda^2(TW).$ (We include the negative sign because if we associate to $\tilde{\Omega}_{\eps}(m)$ a nondegenerate skew-symmetric matrix, then to the dual bivector at $m$ we associate the negative inverse of this matrix.) These bivectors are nondegenerate Poisson structures and we may rewrite equation (\ref{eq:dyn}) as
\begin{equation}
\labell{eq:dyn2}
X_{\eps}= - \tilde{\Omega}_{\eps}^{-1} [d(\epsilon^{-2}\Phi^* H)].
\end{equation} 

In contrast to the $\tilde{\Omega}_{\eps}$ the Poisson structures $\tilde{\Omega}_{\eps}^{-1}$ do have a well defined limit, $({\Omega}^N)^{-1}$. This is a degenerate Poisson structure whose symplectic leaves are the fibres of $N$. As well, for $\eps \ra 0$ we have $\epsilon^{-2}\Phi^* H \ra  d^2_N H$. Hence, 
\begin{eqnarray*}
X_0 &:=& \lim_{\eps \ra 0} X_{\eps}\\ 
    &=& -({\Omega}^N)^{-1}[ (d^N(d^2_N H)],
\end{eqnarray*}
where $d^N$ denotes the exterior derivative with respect to just the fibre variables. This can be rewritten as $$i_{X_0}{\Omega}^N=d^N(d^2_N H).$$ Indeed, this equation defines the limiting vector field $X_0$ globally and we note that the convergence of $X_{\eps}$ to $X_0$ is $C^k$ for any $k$.

By equation (\ref{eq:quasi}) the flow of $X_0$ is fibrewise quasiperiodic and we have at least $(n-l)$ periodic orbits in each fibre.
\end{proof}

By Lemma \ref{Lemma:lemma1}, the flow on the level $\{H=\eps^2\}$ can be viewed
(up to parameterization) as a small (Hamiltonian) perturbation of the 
quasiperiodic flow of $X_0$ on $\{d^2_N H=1\}$. Hence, one may expect that
under this perturbation the set of periodic orbits of $X_0$ on $\{d^2_N H=1\}$ (called the normal modes of $X_0$ following \cite{we3}) 
splits into periodic orbits on $\{H=\eps^2\}$ whose number is bounded from below
by the cup-length or the sum of Betti numbers of this set. These invariants depend on the pair of fibrewise forms $d^2_N H$ and $\Omega^N$ but should be greater than
or equal to the corresponding invariants of $M$.

When the eigenvalues $a_i(x)$ do not bifurcate as functions
of the parameter $x\in M$, the set of normal modes of $X_0$ is an orbifold, 
\cite{ke}. In this case the perturbative analysis can indeed be carried out
by adapting Moser's method, see \cite{bo,mo2}. This leads to a lower bound on
the number of periodic orbits in terms of the cup-length, \cite{ke}. Furthermore,
when for every $x$ the eigenvalues $a_i(x)$ are equal to each other, a lower
bound in terms of the sum of Betti numbers of $M$ has also been obtained in
\cite{gk1}. (This is the case, for example, when $\codim M=2$. In
particular, it is true for symplectic magnetic flows on surfaces, \cite{gi1,gi4}. The condition is also satisfied for symplectic magnetic flows in higher dimensions when $g=\Omega(\cdot, J \cdot)$ for some almost complex structure $J$ on $M$, \cite{ke}.)

In general, this perturbative approach encounters serious difficulties arising from
the fact that the set of normal modes of $X_0$ may fail to be a manifold
or an orbifold. In this work the limiting dynamics is used only as motivation.

\section{Simplification of the variational problem}

Expecting some of the normal modes of the limiting vector field $X_0$ to always persist under the Hamiltonian perturbation $X_{\eps}$, we will restrict our search for periodic orbits to small loops near $M$. This will considerably simplify the original variational problem.

To begin with, we fix some geometric structure on $TW|_M$. Let $J$ be  an almost complex structure that is compatible with $\Omega$. This yields the Riemannian metric $g_J=\Omega(\cdot, J \cdot )$ on $W$. With the splitting 
$$
TW|_M=TM \oplus (TM)^{\Omega},
$$ 
we then have the decomposition 
$$ 
(T_mW, \Omega, J, g_J) = (T_m M \oplus (T_mM)^{\Omega} , \omega_T \oplus \omega_N, J_T \oplus J_N, g_T \oplus g_N)
$$ 
where the subscripts $T$ and $N$ denote tangential and normal components, respectively. With respect to this splitting we will write $z=(x,y)$ for 
$z \in T_mW$, where $x \in T_mM$ and $y \in (T_mM)^{\Omega}.$

\subsection{A Darboux family}

Following \cite{we4}, we note the existence of a Darboux family for $M \subset W$. This a parameterized version of a Darboux chart. It consists of a neighborhood $U$ of the zero section in $TW|_M$ and a mapping 
$$
\Phi \colon U  \longrightarrow  W
$$ 
onto a neighborhood $V$ of $M\subset W$ such that the following conditions hold.

\begin{enumerate}
\item $U_m = U \cap T_mW$ contains the origin.
\item $\Phi_m = \Phi|_{U_m}$ is a symplectomorphism from $(U_m,\Omega(m))$ to $(V_m, \Omega)$, where $V_m$ is an open neighborhood of $m \in W$.
\item $\Phi_m(0)=m$ and $D_0 \Phi_m $ is the identity. 
\end{enumerate}
 
In addition, we may assume that all neighborhoods $U_m$ are open balls of a fixed radius with respect to  $g_J$. For sufficiently small $\eps > 0$, we may also assume that the level set $\{H=\eps^2\}$ lies in $V$.

\subsection{A new loop space}

Denote by $C_0^{\infty}(S^1,V)$ the open subset of $C^{\infty}(S^1,W)$ consisting of small loops contained in $V$. This is an open neighborhood of the constant loops in $V$. Based on our analysis of the limiting dynamics, this is also where we expect to find low energy periodic orbits.

Now, let ${\Lambda}_m$ be the space of $C^{\infty}$ loops in  $T_m W$ whose projections to $T_m M$ have zero mean, and consider the Fr\'{e}chet space bundle  
$$
{\Lambda}=\bigcup_{m \in M} {\Lambda}_m.
$$ 
The map $\Phi$ pulls back $C_0^{\infty}(S^1,V)$ onto an open neighborhood 
${\mathcal U}$ of the zero section in ${\Lambda}$, which in each fibre ${\Lambda}_m$ consists of loops contained in $U_m$. This follows from the inverse function theorem and, put another way, is  essentially the fact that any small loop in $M$ has a unique mean value in $M$ with respect to the map $\Phi$. To be more precise, for every $\sigma\in C_0^\infty(S^1,V)$, there exists
a unique $m\in M$ such that $\sigma=\Phi_m(z)$ for some loop $z\in\Lambda_m$
which takes values in $U_m$.

 Since $\Omega$ is exact in a neighborhood of any of the loops in $C_0^{\infty}(S^1,V)$, 
the differential form $\frak{F}|_{C_0^{\infty}(S^1,V)}$ is exact and so is 
its pullback by $\Phi$. (Indeed, the second term in \eqref{eq:frakF} is always
exact and  the symplectic area of a small disc bounded by the loop can be
taken as a primitive of the first term.) 
We denote the primitive of $\Phi^*(\frak{F}|_{C_0^{\infty}(S^1,V)})$ by $F^0$ and remark that 
$F^0(m,z)= F^0_m(z)$, where $F_m^0$ is the restriction to ${\mathcal U} \cap {\Lambda}_m$ and is given by 
$$
F^0_m (z) = \int_0^1 \haf g_J(m) ( -J(m) \dot{z}, z)  \, dt - 
\int_0^1 H(\Phi_m (z(t))) \, dt.
$$
This is just the standard action functional for the Hamiltonian
$\Phi_m^*H$ which is defined on the open subset $U_m$ of the symplectic vector space $(T_mW, \Omega (m))$. In particular, as above, the first term is the symplectic area of the disc bounded by $z$ in this space.

We have thus simplified the original variational problem to that of finding a critical point of the functional $F^0$ in ${\mathcal U}$ (a neighborhood of the zero section in the Fr\'{e}chet space bundle $\Lambda$). 

\section{An outline of the proof}

\subsection{Step 1}First we utilize the freedom to choose another Hamiltonian $\tilde{H}\colon W \ra \R$ that shares the level set $\{H=\eps^2\}$ with $H$. In fact, we make our changes to the pullbacks $\Phi_m^*H$ in such a way that each new function $h_m$ keeps the level set $\{\Phi_m^* H=\eps^2\} \subset T_m W$ and is equal to the pull back, by $\Phi_m$, of a new global Hamiltonian $\tilde{H}$ defined on $V \subset W$. We then consider the functional $F(m,z)= F_m(z)$ defined on ${\mathcal U} \subset {\Lambda} $ by  
\begin{equation}
\labell{eq:func}
F_m (z) = \int_0^1 \haf g_J(m) ( -J(m) \dot{z}, z)  \, dt - \int_0^1 h_m(z(t)) \, dt.
\end{equation}

Choosing extensions of the maps $h_m$ so that the functional $F$ is defined on all of ${\Lambda}$ we then prove
\begin{Lemma}
\labell{lemma:1}
The choices and extensions above can be made in such a way that any critical point $z_c$ of $F$ on $\Lambda$, satisfying $F(z_c)>0$, corresponds to a periodic orbit of the original system on $\{H=\eps^2+\rho \frac{\eps^2}{4}\}$ for some $\rho \in [-\eps,\eps]$.
\end{Lemma}

\subsection{Step 2}We extend the domain of definition of $F$ to be the Hilbert space bundle 
$$
\tilde{\Lambda}=\bigcup_{m \in M} \tilde{\Lambda}_m,
$$ 
where $\tilde{\Lambda}_m$ is the space of $H^{\haf}$ loops in $T_mW$ whose projections to $T_m M$ have zero mean.

Theorem \ref{thm:main} will then follow from

\begin{Lemma}
\labell{lemma:2}
There exists a critical point $z_c \in {\Lambda} \subset \tilde{\Lambda}$ of $F$ such that $F(z_c)>0$.
\end{Lemma}

\section{Step 1}

\subsection{ A new Hamiltonian and functional}
In looking for periodic orbits on the level set $\{H=\eps^2\}$ we may replace $H$ by any other function which shares this level set. Here we construct such a function which also shares with $H$ all the level sets close to $\{H=\eps^2\}$. This is accomplished by altering (and extending) the pullbacks $\Phi_m^*H$ into a family of functions $h_m \colon T_mW \ra \R^+$ which is smooth in $m$ and satisfies the following conditions.
\begin{enumerate}

\item There is a function $\tilde{H} \colon V \subset W \ra \R$ such that $\Phi_m^*\tilde{H}=h_m|_{U_m}$ for all $m \in M$ and $\tilde{H}$ shares the level sets $\{ H = \eps^2 + \rho \frac{\eps^2}{4}\}$ with $H$, for all $\rho \in [-\eps ,\,\eps]$.
\item All the $h_m =0$ on an open neighborhood of $T_mM \times \{0\} \subset T_mW$.

\item Let $Q_m(z) = \frac{q}{2}\|y\|_m^2$ for some positive $q$ to be specified
later.
Then $\|\nabla h_m(z) - \nabla Q_m(z)\|_m$ is bounded and 
$h_m = Q_m$ for large $\|y\|_m$.
\item For the functional $F$ defined in \eqref{eq:func}, a critical point $z_c$ with $F(z_c)>0$ corresponds to a periodic orbit of $X_H$ on $\{H=\eps^2+\rho \frac{\eps^2}{4}\}$ for some $\rho \in [-\eps,\eps]$.

\end{enumerate}

\begin{Remark}
\label{rmk:q1}
The constant $q>0$ is chosen so that $q$ is not an even integer and $q$ is 
greater than a certain constant depending only on $W$ and $M$.
The assumption that $q$ is not an even integer is crucial in verifying the 
Palais--Smale condition for the functional $F$ (Claim \ref{claim:ps}). The lower bound
for $q$ is essential in Proposition \ref{prop:linking}; see also
Remark \ref{rmk:q}.
\end{Remark}

\subsubsection{Construction of the $h_m$.}
For $\rho \in [-1,1]$ and for each $m \in M$ let $$S_{\rho,m} =\{ \phi_m^* H=\eps^2 + \rho \frac{\eps^2}{4}\} \subset U_m.$$ Extend the hypersurfaces $S_{\rho,m}$ outside $U_m$ by smoothly and quickly joining them to the hypersurfaces 
$$
\{(x,y)\in T_mW \,\mid \, \|y\|_m = \eps^2 + \rho \frac{\eps^2}{4} \}.
$$ 
We still refer to these extended hypersurfaces as the $S_{\rho,m}$ and note that they will be level surfaces of our new functions $h_m$. Denote the union of these hypersurfaces in $T_mW$ by $$C_m= \bigcup_{\rho \in [-1,1]} S_{\rho,m}.$$ Set  
$$
\gamma_m = \max_{(x,y) \in S_{1,m}} \|y\|_m
\quad\text{and}\quad
 \gamma = \max_{m \in M} \gamma_m.
$$ 

Fixing $q$ as in Remark \ref{rmk:q1}, we choose constants $r$ and $b$ such that 
$$
\gamma <r < 2\gamma 
\quad\text{and}\quad
\frac{q}{2}\pi r^2 <b < q\pi r^2.
$$

We then use the following smooth functions to specify the behavior of $h_m$ in $C_m$ and asymptotically in the normal directions. Let $f \in C^{\infty}( [-1,1], \R^+)$ have the properties
$$
f(s)=\begin{cases}
0 &\text{for $s\in (-1,-\eps]$}\\
b & \text{for $s\in [ {\eps}, 1)$}=b
\end{cases}
$$
and 
$$
f'(s)>0 \quad \text{for}\quad -{\eps} <s< {\eps}.
$$
Also, let $g \in C^{\infty}( (0, \infty), \R^+)$ satisfy
$$
g(s)=
\begin{cases}
b & \text{for $s \leq r$}\\
\frac{q}{2}\pi s^2 & \text{for large $s$}
\end{cases}
$$
and 
$$
g(s) \geq \frac{q}{2}\pi s^2\quad\text{and}\quad
0 < \dot{g}(s) \leq q \pi s \quad\text{for}\quad s>r.
$$
Note that $T_m W \backslash C_m$ has two connected components, $A_m$ and $B_m$, where $B_m$ is the open set containing $T_m M \times \{0\} \subset T_m W$. Finally, we set 
$$
h_m(z)=
\begin{cases}
0 & \text{if $z \in  B_m$}\\
f(\rho) & \text{if  $z\in S_{\rho,m}$ for  $-1<\rho<1$}\\
b & \text{if $z \in A_m$ and  $\|y\|_m \leq r$}\\
g(\|y\|_m) &\text{if  $\|y\|_m \geq r$}.
\end{cases}
$$
(See Figure 4.1 below). The functions $h_m$ defined in this way clearly satisfy conditions $1,\,2$ and $3$ as stated above.

\begin{figure}[h]
  \caption{The functions $h_m$ }
  \center{\epsfig{file=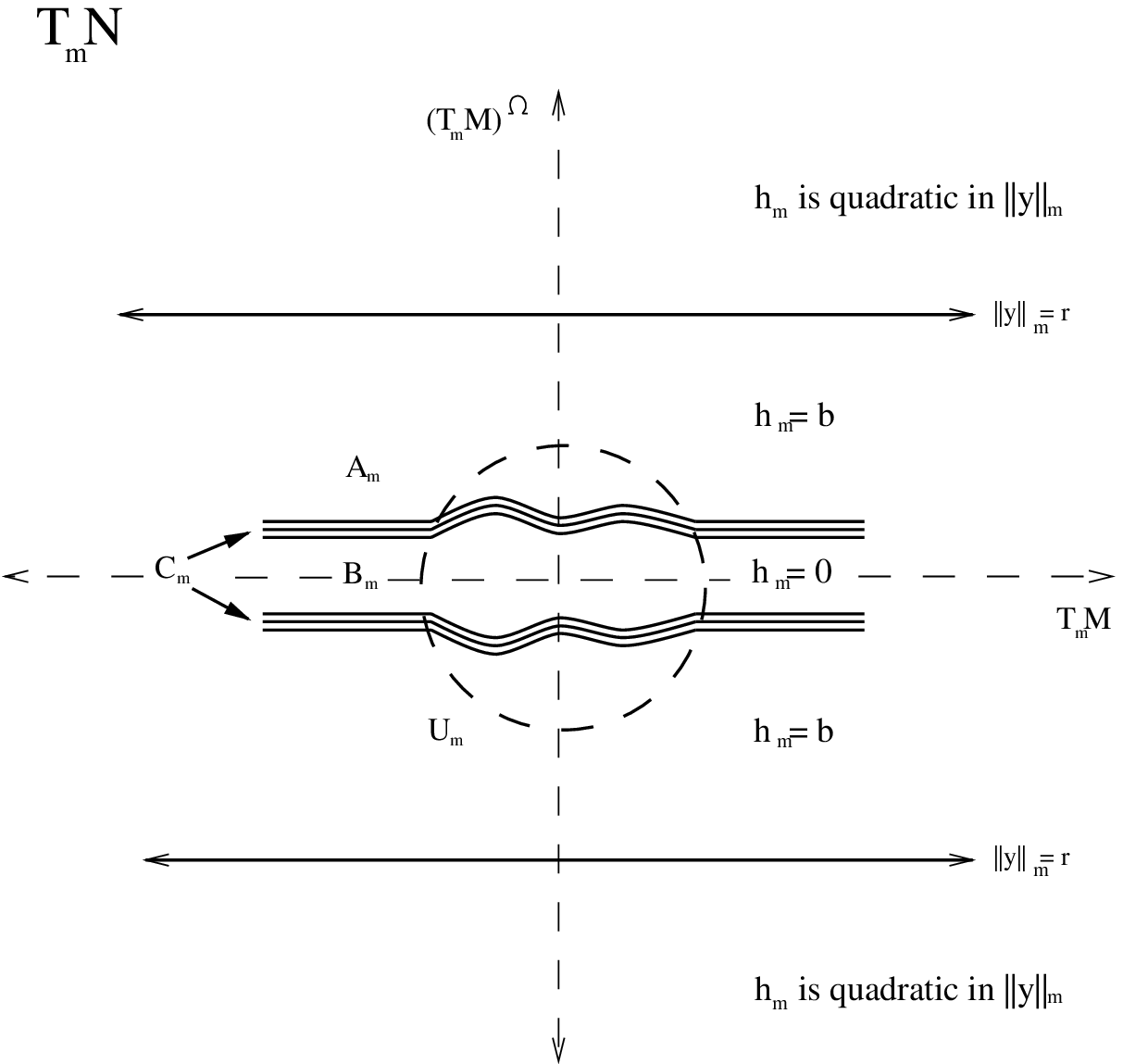}}
\end{figure}

The choices of $f$ and $g$ also yield the following inequalities which will be used later:
\begin{equation}
\labell{eq:u1}
-b + \frac{q}{2} \pi \|y\|_m^2 \leq h_m (z) \leq \frac{q}{2} \pi \|y\|_m^2 + b,
\end{equation}
\begin{equation}
\labell{eq:u2}
\|\nabla h_m (z)\|_m \leq c_1 \|z\|_m \quad\text{for some}\quad c_1 \in \R^+.
\end{equation} 
The second of these is the most crucial and follows from the inequality 
$$
\|\nabla_N h_m (z)\|_m \leq q \pi \|y\|_m \quad\text{for}\quad \|y\|_m \geq r
$$ 
coupled with the facts that $h_m (0)=0$ and that $\nabla_T h_m$ has compact support.

\begin{Remark}
Each function $h_m$ is of the type considered in \cite{gi4}. They differ from those in \cite{vi} and \cite{hz1} in that the level sets $S_{\rho,m}$ are not compact and the functions are asymptotically quadratic only in the normal directions. 
\end{Remark}

\subsection{Locating critical points with positive action}

We now prove that the functions $h_m$ have property $4$.

\begin{Claim}
A critical point of $F_m$, say $z(t)=x(t)+y(t)$, satisfying $F_m(z)>0$ must lie in $S_{\rho,m} \cap U_m$ for some $\rho \in [-\eps, \eps]$. 
\end{Claim}

\begin{proof}
First we show that $F_m(z)>0$ implies that $z \in S_{\rho,m}$ for some $\rho \in [-\eps, \eps]$. If $z(t)$ is a constant solution then $$F_m(z)=-\int_0^1 h_m(z(t)) \, dt \leq 0,$$ since $h_m \geq 0$. Hence, we only have to discount those critical points with $\|y(t)\|_m \geq r \text{ for some } t \in [0,1].$ In fact, because $h_m$ depends only on the fibre variable past $r$, we get $\|y(t)\|_m=\|y(0)\|_m$ for all $t$ in $[0,1].$ When $\|y(t)\|_m \geq r$, we also have $$\nabla h_m(z(t))= \dot{g}(\|y(0)\|_m)\frac{y(t)}{\|y(0)\|_m}.$$ Hence,
\begin{align*}
F_m(z) &= \int_0^1 \left( \haf {g}_J (m)(-{J}(m) \dot{z}(t) ,z(t) ) -h_m(z(t)) \right) \,dt\\
       &= \int_0^1 \left(\haf \dot{g} (\|y(0)\|_m)\|y(0)\|_m -h_m(y(0)) \right) \,dt \\
       &= \haf \dot{g} (\|y(0)\|_m)\|y(0)\|_m -h_m(y(0))\\
       &\leq \haf \dot{g} (\|y(0)\|_m)\|y(0)\|_m -\frac{q}{2} \pi \|y(0)\|_m^2\\
       &\leq 0.
\end{align*} 
The last two inequalities follow from the fact that 
$$
h_m(x,y)=g(\|y\|_m) \geq \frac{q}{2} \pi \|y\|_m^2
$$ 
for $\|y\|_m \geq r$, and $\dot{g}(s) \leq q \pi s$.

Next we show that $x(t) \subset U_m \cap T_mM$. A critical point of $F_m$ is a one-periodic solution of the Hamiltonian dynamical system defined on $T_mW$ by $h_m$ and $\Omega(m)$. Thus, $x(t)$ has period one and because of our splitting of $TW|_M$ it satisfies 
$$
\dot{x}(t) = {J}_T(m) \nabla_T h_m(z(t)).
$$ 
Now, since $x(t)$ has zero mean,
$$
x(t) = \int_0^1 (x(s)-x(t)) \, ds.
$$
Hence,
\begin{eqnarray*}
\|x(t)\|_m &\leq& \int_0^1 \|x(s)-x(t)\|_m \, ds\\
&\leq& \sup_{s\in [0,1]}\|\dot{x}(s)\|_m \\
 &\leq&  \|{J}_T(m)\|_m \sup_{z \in T_mW} \|\nabla_T h_m(z)\|_m. 
\end{eqnarray*}
It is easily checked for our choice of $f$ that $\|\nabla_T h_m(z)\|_m$ is of order $\eps$ for all $z \in C_m $. Hence, $x(t)$ remains sufficiently close to the origin and the proof of the claim is complete.
\end{proof}

Since a critical point of $F$ on $\Lambda$ must also be a critical point of $F_m$ on $\Lambda_m$ for some $m \in M$, the claim yields Lemma \ref{lemma:1}.

\begin{Remark}

At this point we can see why a stronger existence result is not attainable using these techniques. For example, consider what happens when we try to prove that there are periodic orbits on level sets arbitrarily close to a fixed level set $\{H=\eps^2\}$. As in the previous claim, we would like to distinguish any periodic orbits on the candidate level sets by forcing them to have positive action. Accordingly, we define the functions $h_m$ using a new function $f$ which switches from $0$ to $b$ on an arbitrarily small neighborhood of zero, say $(-\delta, \delta)$. Unfortunately, when $\delta$ is too small, i.e. when we look for periodic orbits on level sets too close to $\{H=\eps^2\}$, we lose control of the size of $\|\nabla_T h_m(z)\|_m$. Consequently, given a critical loop $z \in \Lambda_m$ with positve action, we no longer know that $z$ lies within $U_m$. Some part of it may lie in the extended portions of $S_{\rho,m}$ and so $z$ no longer corresponds to a periodic orbit of our original system.    

\end{Remark}

\section{Step 2}

\subsection{Extending the domain of F}

We extend the domain of the functional $F$ from the Fr\'{e}chet space bundle $\Lambda$ to the Hilbert space bundle $\tilde{\Lambda}$ defined below. This extension is motivated by the fact that $F$ has a simple form on $\tilde{\Lambda}$ which allows us to easily verify that the negative gradient flow of $F$ has the properties necessary to employ Minimax techniques to detect critical points. In particular, we are able to extend the ``linking argument'' of Hofer and Zehnder for $\R^{2n}$ (see \cite{hz1,hz2,hz3}) to our bundle $TW|_M$ over $M$.
 
\begin{Remark}Alternatively, one can show that $F$ has a positive critical value
by using the cohomological argument as in \cite{vi} or \cite{gi4}, 
combined with the reduction to finite dimensions from \cite{cz}.
\end{Remark} 

Let $\tilde{\Lambda}$ be the Hilbert bundle over $M$ with fibres $\tilde{\Lambda}_m$ consisting of $H^{\haf}$ loops in $T_mW$ whose projections to $T_mM$ have zero mean. We may consider $\tilde{\Lambda}_m$ as the space of Fourier series $$z(t)=\sum_{k \in \Z} e^{k2\pi J(m) t} z_k,$$ with $z_k=x_k+y_k \in T_m W$ and $x_0=0$, which converge with respect to the norm $\|\,\,\,\|_{\haf,m}$ given by the inner product $$ \l z,z' \r_m =  g_J (m)(z_0,{z'}_0) + 2\pi \sum_{k \in \Z} |k|g_J (m)(z_k,z_k').$$This space clearly includes ${\Lambda}$ and is contained in the bundle of fibrewise $L^2$-loops.

The bundle $\tilde{\Lambda}$ splits naturally in two ways. First there is the orthogonal splitting which in each fibre has the form 
$$
\tilde{\Lambda}_m = E_m^- \oplus E_m^0 \oplus E_m^+ .
$$ 
Here, the space $E_m^-$ consists of the series with nonzero Fourier coefficients for $k<0$ only. The spaces $E_m^0$ and $E_m^+$ are defined analogously. We also have the orthogonal splitting of the $\tilde{\Lambda}_m$ into loops contained in $T_m M$ and loops contained in $(T_mM)^{\Omega}$. For example, $z(t) = x(t) + y(t)$ where the Fourier coefficients of $x(t)$ and $y(t)$ are contained in $T_m M$ and $(T_mM)^{\Omega}$, respectively. We denote this splitting by  $$\tilde{\Lambda} = E_T \oplus E_N.$$

In considering the functional $F$ on the bundle $\tilde{\Lambda}$ we focus first on the fibres, where we have $$F_m (z) = \int_0^1 \haf \tilde{g}_J(m) ( -J(m) \dot{z}, z)  \, dt - \int_0^1 h_m (z(t)) \, dt.$$ With respect to the orthogonal splitting $z=z^- + z^0 + z^+$ it is straightforward to check that 
$$ 
F_m(z) = \frac{1}{2}( \|z^+\|_{\haf,m}^2 -\|z^-\|_{\haf,m}^2) - \int_0^1 h_m (z(t)) \, dt 
$$ 
and the $H^\haf$-gradient with respect to the fibre variables of $\tilde{\Lambda}_m$ is 
$$
\nabla F_m(z)(v) = g_J(m)( z^+-z^-,v) -\int_0^1 {g}_J (m) (\nabla h_m (z) ,v ) \, dt.
$$  

The total gradient flow of $F$ is actually comprised of the gradient flows of the 
$F_m$ on the fibres $\tilde{\Lambda}_m$, coupled with a smooth flow on $M.$ Since $M$ is 
compact, the behavior of the flows on the fibres is the only essential
component in considering compactness properties of the total flow. 

\begin{Claim}
\labell{claim:global}
The vector field $\nabla F$ is smooth and has a globally defined flow on $\tilde{\Lambda}$.
\end{Claim}
The smoothness of $\nabla F$ follows from that of the functions $h_m$. 
Inequality (\ref{eq:u2}) then implies that $\nabla F$ is sublinear in the fibre directions which yields the completeness of the gradient flow (see \cite{amr}).

\begin{Claim}
\labell{claim:smooth}
Critical points of $F_m$ on $\tilde{\Lambda}_m$ are smooth.
\end{Claim}
This is a standard regularity result (see \cite[Lemma 5, p. 88]{hz3}). It justifies the extension of domains since any critical point of $F$ must satisfy $\nabla F_m (z) =0$ for some $m \in M$ and so $\nabla F(z)=0$ implies that $ z \in \Lambda.$ 

\begin{Claim}
\labell{claim:ps}
$F$ satisfies the Palais-Smale condition on $\tilde{\Lambda}$, provided that
$q$ is not an even integer.
\end{Claim}

For the sake of clarity we defer the proof of this claim to an appendix. We just mention here that the claim follows from the careful choice of the asymptotic quadratic behavior of the functions $h_m$. 

\subsection{ Proof of Lemma \ref{lemma:2} }

We use the Minimax Lemma to establish the existence of a critical point of $F$ with a positive critical value and so we recall the setting of this theory. Let $G$ be a $C^1$ function on a Hilbert manifold $L$ and let 
$\frak{T}$ be a family of subsets $T \subset L$. Set 
$$
c(G,\frak{T})=
\underset{T \in \frak{T} \, z \in T}{\inf\,\sup} \, G(z).
$$ 
We then have the following (see, e.g. \cite{hz3}):

\begin{Lemma}[Minimax Lemma]
Let the following properties hold for $G$ and $\frak{T}$:
\begin{enumerate}
\item G satisfies the Palais-Smale condition.
\item The gradient vector field of $G$ gives rise to a global flow 
$\psi^t$.
\item The family $\frak{T}$ is positively invariant under the gradient flow, i.e., $\psi^t(T) \in \frak{T}$ for all $T \in \frak{T}$ and $t \geq 0$.
\item $-\infty < c(G,\frak{T}) <\infty $.
\end{enumerate}
Then there exists $z_c \in L$ such that $$\nabla G (z_c)=0 \text{   and   }G (z_c) = c(h,\frak{T}).$$
\end{Lemma}

Denoting the flow of the negative gradient field of $F$ by $\psi^t$ we must now define a $\psi^t$--invariant family of sets $\frak{T}$ such that $0< c(F,\frak{T}) <\infty.$ To achieve this, we extend the linking argument of \cite{hz1}.

\begin{Proposition}
\label{prop:linking}
Let $e^+_N$ be a nonvanishing section of smooth loops in $E_{N}^+$ with $\|e^+_N(m)\|^2_{\haf,m}= 2\pi$ for all $m \in M$ and let $q$ be greater than
a certain constant depending on $e^+_N$ only. Then
there exists a sufficiently large $\tau >0$ such that for all $m \in M$ the 
subsets 
$$\Sigma_m = \{ x^-+y^-+y^0+se^+_N(m)  \in \tilde{\Lambda}_m\, \mid \, \|x^-+y^- + y^0 \|_{\haf,m} \leq \tau,\, 0 \leq s \leq \tau \}$$ satisfy ${F|}_{\partial \Sigma_m} \leq 0$. 

\end{Proposition}

\begin{proof}

First we label the parts of $\partial \Sigma_m$ as follows
\begin{align*}
&\sigma_1 = \{s=0\},\\
&\sigma_2 = \{s=\tau\},\\
&\sigma_3 = \{ \|x^- + y^- + y^0 \|_{\haf,m} = \tau \}.
\end{align*}
For $z=x^- + x^+ + y^- + y^0 + y^+ \in \tilde{\Lambda}_m$ recall that
\begin{align*}
F_m(z) &= \haf(\|x^+ + y^+ \|_{\haf,m}^2 - \|x^- + y^- \|_{\haf,m}^2 ) - \int_0^1 h_m(z(t)) \, dt\\
       &= \haf(\|x^+\|_{\haf,m}^2 + \|y^+ \|_{\haf,m}^2) - \haf(\|x^-\|_{\haf,m}^2 + \|y^- \|_{\haf,m}^2)  - \int_0^1 h_m(z(t)) \, dt.
\end{align*} 
On $\sigma_1$ we have points of the form $z=x^- + y^- + y^0 $ so that $$ {F_m|}_{\sigma_1}(z)=  - \haf(\|x^-\|_{\haf,m}^2 + \|y^- \|_{\haf,m}^2)  - \int_0^1 h_m(z_1(t)) \, dt \leq 0.$$

For the other parts of the boundary we need to employ half of inequality \eqref{eq:u1}, namely $$ h_m(z) \geq \frac{q}{2} \pi \|y\|_m^2 -b.$$ This yields
\begin{align*}
\int_0^1 h_m(z(t)) \, dt &\geq \frac{q}{2} \pi \left( \int_0^1 \|y^-(t)\|^2_m\,dt + \int_0^1 \|y^0\|^2_m\,dt +  \int_0^1 \|y^+(t)\|^2_m\,dt \right)- b
\end{align*}
which when restricted to $\Sigma_m$ becomes
\begin{align*}
\int_0^1 h_m(z(t)) \, dt &\geq \frac{q}{2} \pi \left( \int_0^1 \|y^-(t)\|^2_m\,dt + \|y^0\|^2_m + s^2 \int_0^1 \|e_N^+(m)(t)\|^2_m\,dt  \right)- b.
\end{align*}
Overall, on $\Sigma_m$ we now have $$F_m (z) \leq b - \haf \|x^- + y^- + y^0 \|_{\haf,m}^2 - 
s^2(\frac{q}{2} \pi \int_0^1 \|e_N^+(m)(t)\|^2_m\,dt - \pi).$$

For our nonvanishing section $e_N^+$ it is clear that there exists a real constant $c >0$ such that $\int_0^1 \|e_N^+(m)(t)\|^2_m\,dt \geq c$ for all $m \in M$. We now choose $q$ to be greater than $\frac{2}{c}$. Then for $\tau$ large enough to satisfy both ${\tau}^2(\frac{q}{2} \pi c - \pi) \geq b$ and $\haf \tau^2 \geq b$, we have $F_m |_{\sigma_2, \sigma_3} \leq 0.$
\end{proof}

\begin{Remark}
\label{rmk:q}
One can show that $e^+_N$ can be chosen in such a way that it suffices to take $q$
strictly greater than $2l/(n-l)$.
\end{Remark}

\begin{Proposition}
\labell{prop:gamma}
There exists a sufficiently small $\alpha > 0$ such that for all $m \in M$ the subsets $$\Gamma_m = \{ y^+ \in \tilde{\Lambda}_m\, \mid \,  \|y^+\|_{\haf,m}^2 = \alpha \}$$ satisfy ${F |}_{\Gamma_m} \geq \beta >0$ for some $\beta \in \R^+$.
\end{Proposition}

\begin{proof}
Since the functions $h_m$ are equal to zero on a neighborhood of $T_mM \times \{0\} \subset T_mW$, the functions 
\begin{align*}
\beta_m \colon \tilde{\Lambda}_m &\to \R \\
z &\mapsto \int_0^1 h_m(z(t)) \, dt
\end{align*}
satisfy $ \beta_m(0)=0$, $\beta'_m(0)=0$, and $\beta''_m(0)=0 $ for all $m$ in $M$. (Here `` $^\prime$ '' denotes a fibrewise derivative.) Restricting to the spaces $E^+_N(m)$ we see then that $${F_m|}_{E^+_N(m)}(z) = \haf(\|y^+\|^2_{\haf,m}) + O( \|y^+\|^3_{\haf,m}).$$ 

\end{proof}

For $\tau > \alpha,$ each $\Sigma_m$ and $\Gamma_m$ intersect at $ \sqrt{ \frac{\alpha}{2\pi}} e_N^+(m).$ Since $\Sigma_m$ and $\Gamma_m$ depend smoothly on $m$ we may extend them to form global subsets of the  bundle $\tilde{\Lambda}.$ We denote these subsets as $\Sigma$ and $\Gamma$, and note that they also intersect (in each fibre). Since $F |_{\partial \Sigma} \leq 0$ and $F |_{ \Gamma} > 0$, we expect the image of $\Sigma$ under the negative gradient flow of $F$ to still 
intersect $\Gamma.$

\begin{Proposition}
$\psi^t(\Sigma) \cap \Gamma \neq \emptyset$, for all $t \geq 0$.
\end{Proposition}

\begin{proof}

Let $P^-$, $P^0$ and $P^+$ be the projection maps corresponding to the splitting $\tilde{\Lambda}= E^- \oplus E^0 \oplus E^+$. Consider the maps 
\begin{align*}
\varphi^t \colon \Sigma &\to E^- \oplus E^0 \oplus \R e^+_N \\
(m,z) &\mapsto \left( \psi_m^t (m,z), (P^- + P^0)\psi_z^t(m,z)+ (\|\psi_z^t (m,z)\|_{\haf,m}^2-\alpha)e^+_N (\psi_m^t (m,z)) \right),
\end{align*}
where $z=z^-+z^0+se^+_N (m)\in \Sigma_m$.
Letting $Z$ be the zero section in  $E^- \oplus E^0 \oplus \R e^+_N$ it follows easily that $$\psi^t(\Sigma) \cap \Gamma \neq \emptyset \iff {(\varphi^t)}^{-1}(Z)\neq \emptyset. $$

The maps $\varphi^t$ are smooth and we claim that they are also nonlinear Fredholm maps with index equal to zero. This is easily confirmed for $\varphi^0$ and so must also hold for $\varphi^t$ with $t>0$ sufficiently small. Then, since $\varphi^t= (\varphi^{\frac{t}{k}})^k$ for any $k \in \N$, it also holds for arbitrary $t$.

One may also easily verify that $\varphi^0$ is transversal to $Z$ and that ${(\varphi^0)}^{-1}(Z)= e_N^+ \subset \Sigma \cap \Gamma$.

We now consider a fixed $t >0$. Note that if $\varphi^t$ is not transversal to $Z$, then the proof is complete since transversality must fail at some $(m,0)\in Z$ with ${(\varphi^t)}^{-1}((m,0)) \in \Sigma$. Hence, we may assume that $\varphi^t$ is transversal to $Z$. By the extension theorem of Smale (see \cite[Theorem 3.1]{sm}) we can then perturb the homotopy 
$\varphi^r$, $r \in [0,t]$, from $\varphi^0$ to $\varphi^t$ to a transversal 
Fredholm homotopy
\begin{align*}
\tilde{\varphi} \colon \Sigma \times [0,t] &\to E^- \oplus E^0 \oplus \R e^+_N 
\end{align*}
such that $\tilde{\varphi}(\cdot\,,0) = \varphi^0(\cdot)$ and $\tilde{\varphi}(\cdot\,,t) = \varphi^t(\cdot).$
The sign of $F$ on $\partial \Sigma$ and $\Gamma$ ensures that ${(\varphi^s)}^{-1}(Z) \cap \partial \Sigma = \emptyset$ for all $s \geq 0$. Consequently, ${(\tilde{\varphi})}^{-1}(Z) \cap \{ \partial \Sigma \times [0,\,t] \} = \emptyset$. With this, Theorem 3.3 of \cite{sm} implies that ${(\tilde{\varphi})}^{-1}(Z)$ is a $2l+1$ dimensional submanifold of $\Sigma \times [0,\,t]$ with boundary equal to 
$$
\{ {(\varphi^0)}^{-1}(Z) \times \{ 0 \} \} \coprod \{ {(\varphi^t)}^{-1}(Z) \times \{ t \} \}.
$$ 
Upon projecting to $M$ we see that these boundary components must generate 
the same homology class in $H_{2l}(M;\Z_2)$. It is easy to see that the 
first component generates the fundamental class and hence the second 
component must be nonempty. 
\end{proof}

\begin{Remark}
A similar argument shows that $\psi_t(\Sigma)\cap \Gamma_m \neq\emptyset$ for
any $m\in M$. 
\end{Remark}

We now complete the proof of Lemma \ref{lemma:2} and hence Theorem 
\ref{thm:main} with a direct application of the Minimax Lemma. Consider the 
family of subsets $\frak{T}=\{\psi^t(\Sigma) \, \mid \, t \geq 0 \}$. It is 
clearly positively invariant under the flow. Let 
$$
c(F,\frak{T})=\underset{t \geq 0}{\inf} \,\underset{z \in \psi^t(\Sigma)}{\sup} F(z).
$$ 
By the previous propositions and the fact that $F$ takes bounded sets to 
bounded sets we have 
$$
\beta \leq \underset{z \in \Gamma}{\inf} F(z) \leq 
\underset{z \in \psi^t(\Sigma)}{\sup} F(z) \leq \infty.
$$ 
This implies that 
$$
0 <\beta \leq c(F,\frak{T}) < \infty
$$ 
and by the Minimax Lemma we have proven the existence of the desired critical point.

\section{Appendix: the Palais--Smale condition for F}

\begin{proof}[Proof of Claim \ref{claim:ps}:\text{ $F$ satisfies the Palais-Smale condition on $\tilde{\Lambda}$}]

We will make use of two distance functions on the bundle $\tilde{\Lambda}$, which we now define. Consider the two fibrewise norms, $\|\,\,\|_{0,m}$ and $\|\,\,\|_{\haf,m}$, which are given by 
$$ \| z\|_{\frac{s}{2},m} =  \|z_0\|_m^2 + 2\pi \sum_{k \in \Z} |k|^s \|z_k\|_m^2$$
for $s=0$ and $s=1$, respectively, where $z(t)=\sum_{k \in \Z} e^{k2\pi J(m) t} z_k$. Each of these norms yields a fibrewise metric which, when coupled with the base metric $g_J$, define an $L^2$ and an $H^{\haf}$ metric on $\tilde{\Lambda}$. We are interested in the distance functions corresponding to these metrics which we will denote by $d_{L^2}$ and $d_{H^{\haf}}$. 

Given a sequence $\{(m_i, z_i)\} \subset \tilde{\Lambda}$ such that for the $H^{\haf}$-gradient we have
\begin{equation}
\labell{eq:critlim}
\nabla F (m_i, z_i) \ra 0
\end{equation}
with respect to the $H^{\haf}$ metric,
we need to show that there exists a convergent subsequence with respect to $d_{H^{\haf}}$. (Since $M$ is compact we already know that there is a convergent subsequence $m_i \ra m_0 \in M$.)

In what follows we focus entirely on the fibre component of the 
gradient (equal to
the gradient of the restriction to a fibre), which we again denote by
$\nabla$. The fibre component, $\nabla F$, is a function 
$\tilde{\Lambda}\to \tilde{\Lambda}$, where we identify a tangent space
to the fiber $\tilde{\Lambda}_m$ with the fibre itself.
Clearly, the norm of the fibre component does not exceed the
norm of the gradient and hence \eqref{eq:critlim} still holds for
the fibre components. 

Consider the special form that condition \eqref{eq:critlim} takes on 
$\tilde{\Lambda}$. 
Recall that the functional $F$ is given as 
$F(m,z)=F_m(z)=\alpha_m(z)-\beta_m(z)$, where
$$
\alpha_m (z) = \int_0^1 \haf \tilde{g}_J(m) ( -J_m \dot{z}, z)  \, dt 
\quad\text{and}\quad
\beta_m(z)= \int_0^1 h_m (z(t)) \, dt.
$$  
With respect to the fibrewise orthogonal splitting $z=z^- + z^0 + z^+$ 
(see Section 4.3.2), it is straightforward to check that 
$$ 
\alpha_m(z) = \frac{1}{2}( \|z^+\|_{\haf,m}^2 -\|z^-\|_{\haf,m}^2) 
$$ 
and, since $\nabla$ denotes the fibre component of the gradient,
$$
\nabla \alpha_m(z) = z^+-z^-.
$$ 
As in \cite[Prop. 5 p. 86-7]{hz3}, one can also show that 
$$
\nabla \beta_m(z)=j_m^* \nabla h_m (z),
$$ 
where $j_m^*$ is the formal adjoint of the inclusion 
$j_m \colon H^{\frac{1}{2}}(S^1,T_mW) \ra L^2(S^1,T_mW)$ and is a 
compact map. Indeed, our second inequality \eqref{eq:u2} implies that 
$\nabla h_m$ takes bounded sets in $L^2(S^1,T_mW)$ to bounded sets so that the map
$\nabla\beta\colon \tilde{\Lambda}\to \tilde{\Lambda}$ defined
as $\nabla\beta(m,z) = \nabla \beta_m (z)$, is also compact.

We may now rewrite \eqref{eq:critlim}, in a slightly weakened form, as 
$$
\|(z_i^+ - z_i^-)-(\nabla \beta(m_i,z_i))\|_{\haf, m_i} \to 0.
$$ 

Assume first that the $\|z_i\|_{\haf, m_i}$ are bounded for an infinite
subsequence of points. Without loss of generality we may assume that
$(m_i,z_i)$ is this subsequence. By the compactness of the map 
$\nabla \beta$ we know that $\{(m_i,\nabla \beta(m_i,z_i))\}$ is 
relatively compact. Thus, the sequence $\{(m_i,z_i^+ - z_i^-)\}$ has a 
convergent subsequence. Note that $z_i^+$ and $z_i^-$ are orthogonal to 
each other in $\tilde{\Lambda}_{m_i}$. After passing if necessary to
subsequences, this implies that each of the 
sequences $\{(m_i,z_i^+)\}$ and $\{(m_i,z_i^-)\}$ converges for the
same subsequence of points $(m_i,z_i)$. As before, we may restrict our
attention to this subsequence. Finally, since $\{(m_i,z^0_i)\}$ 
is a bounded sequence in a finite dimensional space it too has a 
convergent subsequence and the proof in this case is finished. 

Looking for a contradiction, we assume that $\|z_i\|_{\haf, m_i}$  
are unbounded for some infinite sequence on which we now focus.
Set 
$$
u_i=\frac{z_i}{\|z_i\|_{\haf, m_i}}
\quad\text{and}\quad
w_i=\frac{\nabla h_{m_i}(z_i)}{\|z_i\|_{\haf, m_i}}.
$$ 
The assumption \eqref{eq:critlim} now takes the form 
$$
\left\| u_i^+ -u_i^- -j_{m_i}^* w_i\right\|_{\haf, m_i} \to 0.
$$

By inequality \eqref{eq:u2}, the 
sequence $\|w_i\|_{L^2,m_i}$ is bounded.  Indeed,
$$
\| w_i\|_{L^2,m_i}\leq 
\frac
{\| \nabla h_{m_i}(z_i)\|_{L^2,m_i}}{{\|z_i\|_{\haf, m_i}}}
\leq \frac{\| \nabla h_{m_i}(z_i)_i\|_{L^2,m_i}}{{\|z_i\|_{L^2, m_i}}}
\leq c_1.
$$
The compactness of the operators $j_{m_i}^*$ then implies that 
$\{(m_i,j_{m_i}^* w_i)\}$ 
is relatively compact in $\tilde{\Lambda}$ with respect to 
$d_{H^{\haf}}$. Hence, the sequence $\{(m_i, u_i^+ -u_i^-)\}$ is
also relatively compact in $\tilde{\Lambda}$ with respect to 
$d_{H^{\haf}}$. Just as above we then get a convergent subsequence 
$$
(m_i,u_i) \ra (m_0,u).
$$ 
Note also that
$$
\| u\|_{\haf,m}=\lim \| u_i\|_{\haf,m_i}=1.
$$

Now 
\begin{eqnarray*}
d_{L^2}\big((m_i,w_i), (m_0,\nabla Q_{m_0}(u))\big) 
& \leq &  d_{L^2}\big((m_i,w_i), (m_i,\nabla Q_{m_i}(u_i))\big) \\
&&\quad + d_{L^2}\big( (m_0,\nabla Q_{m_0}(u)), 
(m_i,\nabla Q_{m_i}(u_i))\big),
\end{eqnarray*}
where $Q_m(z)$ is the quadratic term in the definition of $h_m$. 
The second term on the right hand side goes to zero
as $i\to\infty$ because $\nabla Q$ is continuous.
As for the first term, we have 
$$
d_{L^2}\big((m_i,w_i), (m_i,\nabla Q_{m_i}(u_i))\big)= 
\frac{1}{\|z_i\|_{\haf, m_i}} \|\nabla h_{m_i}(z_i)
-\nabla Q_{m_i}(z_i)\|_{L^2,m_i}.
$$ 
However, by the construction of $h_m$, the difference $h_m-Q_m$ has compact 
support, and hence, as is easy to see, 
$ \| \nabla h_m - \nabla Q_m\|_{L^2}$ is bounded on 
$\tilde{\Lambda}$. Thus the first term also goes to zero and 
$$
(m_i,w_i) \ra (m_0,\nabla Q_{m_0}(u))
$$ 
with respect to $ d_{L^2}$. This means that  
$$
(m_i,j_{m_i}^*w_i) \ra (m_0,j_{m_0}^*\nabla Q_{m_0}(u))
$$ 
with respect to $ d_{H^{\haf}}.$ Accordingly, $u$ satisfies 
$$
u^+-u^- - j_{m_0}^*(\nabla Q(m_0,u)) = 0.
$$

This is equivalent to 
$$
{J}_N({m_0}) \dot{u}_N = q \pi u_N \quad\text{and}\quad u_T=0.
$$ 
However, $u_N \ra {J}_N({m_0}) {\dot{u}}_N$ is a self adjoint map on $E_N(m_0)$ with spectrum $\{2\pi \Z\}$. Since $q$ is not a positive even integer, this forces $u_N=0$ and we get a contradiction to $\|u\|_{\haf,m_0}=1$.
\end{proof}

\end{document}